\documentclass[12pt]{article}
\usepackage{amsmath}
\usepackage{amssymb}
\usepackage{amscd}
\usepackage{epsfig}

\def\no{\noindent}
\numberwithin{equation}{subsection}

\newenvironment{theorem}{\par\bigskip\noindent\refstepcounter{equation}{\bf Theorem \theequation.}\em}{\em\par\bigskip\noindent}

\newenvironment{proposition}{\par\bigskip\noindent\refstepcounter{equation}{\bf Proposition \theequation.}\em}{\em\par\bigskip\noindent}
\newenvironment{lemma}{\par\bigskip\noindent\refstepcounter{equation}{\bf Lemma \theequation.}\em}{\em\par\bigskip\noindent}
\newenvironment{corollary}{\par\bigskip\noindent\refstepcounter{equation}{\bf Corollary \theequation.}\em}{\em\par\bigskip\noindent}
\newenvironment{remark}{\par\bigskip\noindent\refstepcounter{equation}{\bf Remark \theequation.}}{\par\bigskip\noindent}

\newcommand{\R}{\mathbb R}

\newcommand{\Z}{\mathbb Z}

\newcommand{\qed}{{\unskip\nobreak\hfil
        \penalty50\hskip1em\hbox{}\nobreak\hfil
        $\square$\parfillskip=0pt\finalhyphendemerits=0 \par}}

\def\ga{\gamma}
\def\Ga{\Gamma}

\begin{document}

\title{Conjugacy rigidity for nonpositively curved graph manifolds}
\author{Christopher B. Croke\thanks{Supported by NSF grant DMS 99-71749 and a CRDF grant.}}
\date{\today}
\maketitle

\begin{abstract}
We show that the metric of nonpositively curved graph manifolds is determined by its geodesic flow.  More precisely we show that if the geodesic flows of two nonpositively curved graph manifolds are $C^0$ conjugate then the spaces are isometric.
\end{abstract}

\setcounter{section}{1}
\setcounter{subsection}{0}

\subsection{Introduction}

Two complete Riemannian manifolds, $N$ and $N'$, are said to have $C^0$ conjugate geodesic flows if there is a continuous map $F:S^1N\to S^1N'$ such that $F\circ g_t=g_t'\circ F$ for all $t$ where $g_t$ (resp. $g_t'$) represents the geodesic flow on the unit sphere bundle $S^1N$ (resp. $S^1N'$).  The conjugacy is called $C^k$ if $F$ is $C^k$.

Graph manifolds of nonpositive curvature are compact three dimensional manifolds whose geometry has been extensively studied in \cite{Sh}.  They are of interest for a number of reasons, but in particular they are the easiest nontrivial examples of nonpositively curved manifolds that have Ballmann rank 1 but whose fundamental group is not hyperbolic.  Here we show:

\begin{theorem}
\label{main}
If $N$ and $N'$ are compact nonpositively curved graph manifolds whose geodesic flows are $C^0$ conjugate then they must be isometric.
\end{theorem}

The conjugacy rigidity question for compact Riemannian manifolds: ``When does a $C^k$ conjugacy between the geodesic flows of two compact Riemannian manifolds $N$ and $N'$ imply that $N$ and $N'$ are isometric?'' is a much studied but still largely open question.  This question is closely related to other rigidity problems (such as the "boundary rigidity" problem and - in negative curvature - the "marked length spectrum" rigidity problem).  

There are many examples where isometry does not hold (Zoll surfaces or more generally the examples in section 6 of \cite{C-K2} for example) however none of these occur in the nonpositive curvature setting.  The only general conjugacy rigidity result known where there are no curvature assumptions made is that of \cite{C-K2} where one of the metrics is assumed to have a parallel vector field (hence it is locally a Riemannian product where one factor is an interval).  A natural setting for these questions is when at least one of the manifolds (say $N$) is assumed to have nonpositive curvature.  The two dimensional case (when $N$ has nonpositive curvature) is completely taken care of by \cite{C-F-F}, which is based on \cite{O}, and \cite{Cr1}.  In higher dimensions, if $N$ is a locally symmetric space of negative curvature (see \cite{B-C-G}) or a flat manifold (see \cite{Cr2}) and if the conjugacy is $C^1$ (or $C^0$ if $Vol(N)=Vol(N')$) then $N'$ is isometric to $N$.

In this paper we consider the case where both manifolds $N$ and $N'$ have nonpositive curvature and $n\geq 3$ (since the 2 dimensional case was taken care of by \cite{C-F-F}).  In this setting it is known that $N'$ must be isometric to $N$ if $N$ has Ballmann rank $\geq 2$ (see \cite{C-E-K}). 

In the above I restricted my attention only to rigidity theorems for a pair of compact manifolds without boundary.  There are a number of results on compact manifolds with boundary or on finite volume manifolds.  Much work has also been done on the deformation question (i.e. can there be nontrivial differentiable one parameter families $N_t$ with the "same" geodesic flows) (see for example \cite{G-K}, \cite{C-S}, \cite{D-S}, and \cite{S-U}).

In a previous paper, \cite{C-K1}, motivated by a question of Gromov in \cite{Gr}, the author along with Bruce Kleiner considered the space of geodesic rays on the universal cover of such nonpositively curved graph manifolds (as well as other related spaces).  The asymptotic properties of the geodesic rays are encoded in the geodesic boundary (equivalence classes of rays) and the induced action of the fundamental group on that boundary. We were able to give necessary and sufficient conditions on the metrics to guarantee that two such (necessarily homeomorphic) spaces have universal covers with equivalent boundaries (i.e. there is a an equivariant, with respect to the action of the fundamental group, homeomorphism between the boundaries).  In the process we saw that there were examples of non isometric graph manifolds of nonpositive curvature whose universal covers have equivalent boundaries.  It thus became natural to ask if complete knowledge of the geodesic flow (i.e. knowledge up to $C^0$ conjugacy) is enough to completely determine the metric.  (Such a conjugacy, of course, induces an equivalence between boundaries in a natural way.)  Theorem \ref{main} gives a positive answer to that.

Graph manifolds of nonpositive curvature have a very restricted geometric form as studied in \cite{Sh}.  They are the union over totally geodesic boundaries of pieces (the ``geometric Seifert components'') whose universal covers are isometric to a Riemannian product $X\times \R$ where $X$ is a surface of nonpositive curvature with geodesic boundary.  The universal cover of the geometric Seifert components are minimal displacement sets (see section \ref{defs} for a definition) of isometries.  In fact, the minimal displacement set of an isometry of any simply connected space of nonpositive curvature has the form $X\times \R$ where $X$ is a convex subspace.  They thus become amenable to some techniques developed in \cite{C-K2}.  In sections \ref{defs} and \ref{mindispsets} of the paper we thus consider $C^0$ geodesic conjugacies between arbitrary nonpositively curved manifolds of dimension $\geq 3$.   In section \ref{mindispsets} we show that for corresponding minimal displacement sets  the conjugacy induces affine maps between the corresponding spaces $X$ (see Proposition \ref{nearisometry}).

In section \ref{graph} we restrict our attention to graph manifolds and prove  Theorem \ref{main}.  Using the results in section \ref{mindispsets} it is straightforward to show that corresponding Seifert components are isometric.  However, we note that it is easy to glue isometric Seifert components together so as to make non isometric spaces.   The main point of section \ref{graph} is to show that the components must glued together in an isometric manner. 

We are grateful to Bruce Kleiner for helpful comments on earlier versions of this paper.

\subsection{Notation and preliminaries}
\label{defs}

A standard consequence of the flat strip theorem says that if $\nu$ is a geodesic in a simply connected nonpositively curved manifold $M$ then the space $M_\nu =\{x\in M|x$ lies on a geodesic that stays a bounded distance from $\nu\}$ is a convex subset of $M$ which is isometric to $X_\nu\times \R$ where $X_\nu$ is a simply connected manifold of nonpositive curvature with convex boundary.  The unit speed geodesics that stay a bounded distance from $\nu$ are precisely the "vertical" curves $t\to (x_0,t+t_0)\in X_\nu\times \R$ for fixed $(x_0,t_0)$.  We will let $\pi_\nu:M_\nu\to X_\nu$ (resp. $\R_\nu:M_\nu\to \R$) be the projection.

Let $N$ be a Riemannian manifold of dimension $\geq 3$ and of nonpositive curvature and $\pi_N:M\to N$ its universal cover.  Since basepoint projection, $b:S^1N\to N$, induces an isomorphism between $\pi_1(S^1N)$ and $\pi_1(N)\equiv\Ga$ we see that $\pi_{N*}:S^1M\to S^1N$ is also a universal cover.  (For any tangent vector $v$ to any space we will let $b(v)$ represent its base point.)  After appropriate choices of base points we get cocompact actions of $\Ga$ on $M$ and $S^1M$ with $M/\Ga=N$ and $S^1M/\Ga=S^1N$.  For every $\ga\in\Ga$ we let $M_\ga\subset M$ be the set of points of minimal $\ga$-displacement, i.e. the set of $m$ such that $d(m,\ga(m))$ is a minimum.  This is exactly the space $M_\nu$ where $\nu$ is any axis of $\ga$ (i.e. $\ga(\nu)=\nu$) and the axes of $\ga$ are precisely the vertical geodesics in $M\ga$. Thus $M_\ga$ is a convex subset of $M$ which is isometric to $X_\ga\times \R$ where $\ga$ acts by translation on the $\R$ factor and $X_\ga$ is a simply connected manifold of nonpositive curvature with convex boundary.  (We will use $\pi_\ga:M_\ga\to X_\ga$ and $\R_\ga:M_\ga\to \R$) to represent the projections.)  The centralizer, $Z(\ga)\subset \Ga$, of $\ga$ in $\Ga$, will act cocompactly on $M_\ga$ preserving the vertical directions.    

For every $m\in M_\ga$, $\pi_N$ takes the (vertical) geodesic segment $\overline{m\ga(m)}$ to a closed geodesic in the free homotopy class of $\ga$, and all such closed geodesics come this way.

\subsection{Geodesic conjugacies and the $M_\nu$}
\label{mindispsets}

In this section we consider compact nonpositively curved manifolds $N$ and $N'$ of dimension $\geq3$.  We discuss the consequences that a $C^0$ conjugacy of the geodesic flows has on the spaces $M_\nu$ (and in particular on the $M_\ga$).

The isomorphism $F_*:\pi_1(S^1N)\to \pi_1(S^1N')$, induces an isomorphism between $\pi_1(N)$ and $\pi_1(N')$.  We will let $\Ga$ represent the common abstract group.  Thus (after choices of base points compatible with $F$ for $M$ $M'$, $N$, and $N'$) we can lift $F$ to a uniformly continuous conjugacy $\tilde F:S^1M\to S^1M'$ and have
 cocompact actions of $\Ga$ on $M$ and $M'$ whose differentials commute with $\tilde F$.  We will let $D(F)<\infty$ be the maximum over $m\in M$ of the diameters of the sets $\tilde F(S^1_mM)\subset S^1M'$ where $S^1_mM\subset S^1M$ represents the unit vectors with basepoint $m$. 

\begin{lemma}
\label{dist}
For all $v,w \in S^1M$:
$$-2D(F) \le d_M(b(v),b(w))-d_N(b(\tilde F(v)),b(\tilde F(w)))\le 2D(F)$$ 

\end{lemma}
\label{bddist}
\noindent
{\bf Proof:}
For $v,w \in SM$, let $\tau :[0,l] \to M$ be
a minimizing geodesic from $b(v)$ to $b(w)$. Then both
$d_N(b(\tilde F(v)),b(\tilde F(\tau'(0))))$, and
$d_N(b(\tilde F(w)),b(\tilde F(\tau'(l))))$, are bounded by $D(F)$ and
$d_N(b(\tilde F(\tau'(0))),b(\tilde F(\tau'(l))))=d_N(b(\tilde F(\tau'(0))),b(g^l(\tilde F(\tau'(0))))) = l = d_M(b(v),b(w))$.
Hence the inequalities follow from the triangle inequality.
\qed

\bigskip

For any geodesic $\nu$ of $M$, $\tilde F$ will take the unit tangent vector field of $\nu$ to the unit tangent field of a geodesic of $M'$ which by abuse of notation we will call $\tilde F(\nu)$.

For $\nu$ a geodesic in $M$ and $m\in M_\nu$ we let $V(m)$ be the vertical direction at $m$ (i.e. points in the positive $\R$ direction of $M_\nu=X_\nu\times \R$).  An immediate consequence of the above lemma is that if a geodesic $\nu_1$ stays a bounded distance from $\nu$ then $\tilde F(\nu_1)$ stays a bounded distance from $\tilde F(\nu)$.   Thus $\tilde F(V(m))=V'(m')$ for some $m'\in M'_{\tilde F(\nu)}$.  This allows us to define a map $G_\nu:M_\nu\to M'_{\tilde F(\nu)}$ with $\tilde F(V(m))=V'(G_\nu (m))$, and then, since we have $G_\nu(\tau_t m)=
\tau'_t(G_\nu(m))$ (where $\tau_t$ and $\tau'_t$ represent vertical translation in $M_\nu=X_\ga\times \R$ and $M'_\nu$ respectively), we get a well defined function $GX_\nu:X_\nu\to X'_{\tilde F(\nu)}$.

For $\ga\in \Ga$ the above applies to $M_\ga$.  Since $F$ must take closed geodesics in $N$ in the free homotopy class of $\ga$ to corresponding closed geodesics (of the same length) in $N'$ we see that $\tilde F$ takes $M_\ga$ to $M'_\ga$ and gives rise to maps $G_\ga$ and $GX_\ga$ as above.  This also tells us that the translation by $\ga$ on $M_\ga$ and $M'_\ga$ is by the same amount.  We also note that since $\tilde F$ commutes with the differentials of $Z(\ga)$, and since $Z(\ga)$ preserves $M_\ga$ and vertical directions, we have that $G_\ga$ commutes with $Z(\ga)$ (which acts cocompactly). 

A curve $c$ in $M_\nu=X_\nu\times \R$ is said to have $\R$ growth rate $A$ if $\lim_{t\to \infty}(\frac {\R_\nu(c(t))} t)=A$.

\begin{lemma}
\label{rrate}
If a curve $c$ in $M_\nu$ has growth rate $A$ then $G_\nu(c)$ also has $\R$ growth rate $A$.
\end{lemma}

{\bf Proof:} 
Let $B_\nu$ be the Busemann function of $\nu$ (i.e. $B_\nu(x)=lim_{t\to \infty}(t-d(x,\nu(t))$).
Now Lemma \ref{bddist} (with both $v$ and $w$ vertical vectors) implies that $|B_\nu(x)-B_{\tilde F(\nu)}(G_\nu(x))|\leq 2D(F)$.  On the other hand $B_\nu(x)=\R_\nu(x)+const$ and hence $\R$ growth rates are the same as Bussemann growth rates (also in $M'_{\tilde F(\nu)}$) and the lemma follows. 

\qed

\begin{proposition}
\label{nearisometry}
For any geodesic $\nu$ in $M$, $GX_\nu$ takes each unit speed geodesic segment in $X_\nu$ to a constant speed geodesic segment in $X'_{\tilde F(\nu)}$.  Further, if the segment is a subset of a ray in $X_\nu$ then the image is unit speed.
\end{proposition}

{\bf Proof:} 

Let $\sigma(t)$ for $0\leq t\leq L$ be a unit speed geodesic segment in $X_\nu$ and for each $i$ let $\sigma_i(t)=(\sigma(\frac {t} {\sqrt{1+i^2}}),\frac {it} {\sqrt{1+i^2}}))\subset M_\nu$ for $0\leq t\leq L\sqrt{1+i^2}$ be a unit speed geodesic in $M_\nu$.  Now $F(\sigma'_i(t))$ will be the tangent vector field, $\tau'(t)$, to a unit speed geodesic $\tau_i(t)$ in $M'$.  Let $p_i\in M'_{\tilde F(\nu)}$ (resp. $q_i\in M'_{\tilde F(\nu)}$) be the unique closest point in $M'_{\tilde F(\nu)}$ to $\tau_i(0)$ (resp. $\tau_i(L\sqrt{1+i^2})$).  The uniform continuity of $F$ guarantees $d(F(\sigma_i'(t)),F(V(\sigma_i(t))))$ 
converges to $0$ ($d$ computed in the tangent bundle) and hence that $\lim_{i\to \infty}d(p_i,\tau_i(0))=0$ and $\lim_{i\to \infty}d(q_i,\tau_i(L\sqrt{1+i^2}))=0$.  
Hence if $\bar \tau_t(t)$ for $0\leq t\leq L\sqrt{1+i^2}$ is the constant speed geodesic from $p_i$ to $q_i$ we see (by the nonpositivity of the curvature) that for every $t$, $d(\tau_i'(t),\bar \tau_i'(t))$, and hence $d(\bar \tau_i'(t),F(V(\sigma_i(t))))$, converges to zero.  In particular for every $t$ we have $d(\pi'_{\tilde F(\nu)}(\bar \tau_i'(t)),GX_\nu(\sigma (\frac {t} {\sqrt{1+i^2}})))$ converges to $0$.
But each $\pi'_{\tilde F(\nu)}(\bar \tau_i'(t))$ is a constant speed geodesic in $X'_{\tilde F(\nu)}$ and hence so is $GX_\nu(\sigma(t))$.

Now let $r;[0,\infty)\to X_\nu$ be a unit speed ray and assume that $GX_\nu (r)$ is a ray with constant speed $a$.  We need to show that $a=1$.   Using the construction above for segments, we first define (by letting $L$ go to $\infty$) for each $i$ the unit speed ray $\tau_i(t)$ whose distance to $M'_{\tilde F(\nu)}$ goes uniformly (since $\tilde F$ is uniform) to 0 as $i$ goes to infinity and hence we can define (by a standard limiting argument) a ray $\bar \tau_i(t) \subset M'_{\tilde F(\nu)}$ such that as $i$ goes to $\infty$ both $d(\bar \tau_i'(t),\tau'_i(t))$ and $d(\pi'_{\tilde F(\nu)}(\bar \tau_i'(t)),GX_\nu(r (\frac {t} {\sqrt{1+i^2}})))$ converge uniformly to 0.  Note that $d(\bar \tau_i(t),G_\nu(\sigma_i(t))$ is uniformly bounded since both curves are at a uniformly bounded distance (depending on $D$) from $\tau_i(t)$, and hence, by Lemma \ref{rrate}, $\bar \tau_i$ has the same $\R$ growth rate as $\sigma_i$ (i.e.$\frac i {\sqrt{1+i^2}}$).  But since $GX_\nu(r)$ has speed $a$ (and hence $GX_\nu(r (\frac {t} {\sqrt{1+i^2}}))$ has speed $\frac a {\sqrt{1+i^2}}$), and $\bar \tau_i$ has speed approaching 1 (the speed of $\tau_i$) we see that $\bar \tau_i$ can have the same $\R$ growth as $\sigma_i$ for all $i$ if and only if $a=1$. 
   \qed  

\bigskip

Maps $GX_\nu$ that satisfy the conclusion of Lemma \ref{nearisometry} are called affine maps.  We now recall some of the standard properties of such maps (e.g. see chapter VI of \cite{K-N}).  For $V\in S^1X_\nu$ let $s(V)$ represent the speed of the $GX_\nu$ image of the geodesic determined by $V$.  Since $s$ is continuous it is locally bounded above and hence $GX_\nu$ is locally Lipschitz (in fact globally in the cocompact case).  Thus $GX_\nu$ is differentiable almost everywhere and at every $x\in X_\nu$ where $GX_\nu$ is differentiable we see $GX_\nu=exp_{GX_\nu(x)}\circ D(GX_\nu)_x \circ exp^{-1}_x$ (this is always true locally, but in our case $exp^{-1}_x$ is globally defined).  In particular $GX_\nu$ is $C^\infty$ and $GX_\nu=exp_{GX_\nu(x)}\circ D(GX_\nu)_x \circ exp^{-1}_x$ holds for all $x$.  This implies that $GX_\nu$ preserves the respective connections and hence that $s$ is invariant under holonomy.  Since $X_\nu$ is convex (with boundary) and simply connected, invariant subspaces of the holonomy action induce deRham splittings as usual.  Thus we get:

\begin{corollary}
For any geodesic $\nu$ in $M$, $GX_\nu$ is $C^\infty$, takes deRham factors to deRham factors, acts homothetically on each non Euclidean deRham factor (and hence is an isometry on each noncompact one), and acts as an affine transformation on the flat factor.  In particular if $X_\nu$ is irreducible, not flat, and noncompact (e.g. in the case of $M_\ga$ when $Z(\ga)$ is not cyclic) then $GX_\nu$ is an isometry.
\end{corollary}

\begin{remark}
It is certainly possible that in the above $GX_\nu$ is always an isometry.  In fact, it is highly likely if the conjugacy rigidity conjecture for nonpositively curved manifolds is true since in this case $X_\nu$ and $X'_{\tilde F(\nu)}$ would be isometric.  
\end{remark} 

\begin{remark}
\label{key}
Consider the case that $N$ is 3 dimensional, and for some $\ga \in \Ga$, $X_\ga$ is noncompact.
Combining lemma \ref{nearisometry} and the corollary  we see that  either $X_\ga$ is a flat strip or $G_\ga$ is an isometry that commutes with $Z(\ga)$ and $M_\ga/Z(\ga)$ is isometric to $M'_\ga/Z(\ga)$.  

Let $O\subset X_\ga$ be a flat strip neighborhood of some boundary geodesic $\tau$, and $\bar O\subset M_\ga$ be $\pi_\ga^{-1}(O)$, which has the metric of a flat plane cross an interval. Let $\alpha\in Z(\ga)$, $\alpha\not=\ga$ be such that $\pi_\ga^{-1}(\tau)\subset M_\ga\cap M_\alpha$ then the geometry of $\bar O$ forces $\bar O\subset M_\ga\cap M_\alpha$. 
\end{remark} 

\subsection{Graph Manifolds}
\label{graph}

We now consider the case where $N=M/\Ga$ is a compact $3$-dimensional graph manifold with a nonpositively curved Riemannian metric with fundamental group $\Ga$ and universal cover $M$.  

By the theorem of 
\cite{Sh}, $N$ has the following structure.  There is a finite connected graph $\cal G$ with vertex set $\cal V$ and oriented edge set $\cal E$ (i.e. each edge shows up twice in $\cal E$).  For each $v\in V$ There is a compact nonpositively curved $3$-manifold $N_v$ with nonempty totally geodesic
boundary (the geometric Seifert components of $N$), and Seifert fibrations 
$p_v:N_v\to X_v$ where
the metric on $N_v$ has local product structure compatible
with the fibration $p_v$, and the $X_v$ are nonpositively 
curved orbifolds; $N$ is the union of the $N_v$.  Each $N_v$ is $M_{\ga_v}/Z(\ga_v)$ for some $\ga_v\in \Ga$.  There is one boundary component $T_e$ of $N_v$ for each $e\in \cal E$ such that $v$ is the initial vertex of $e$.  (Note that a nonoriented edge from $v$ to $v$ will correspond to two boundary components of $N_v$.)  The $N_v$ are glued together as is described by $\cal G$.
That is for each edge $e\in \cal E$, say from $v_1$ to $v_2$, we glue $N_{v_1}$ to $N_{v_2}$ along $T_e$ and $T_{-e}$ respectively. (Note in the case where the metric on $N_{v_1}$ near $T_e$ is a product of an interval with $T_e$ then remark \ref{key} says that the picture looks the same on $N_{v_2}$ near $T_{-e}$ and the gluing identifies the whole product neighborhoods.)  

We note that these gluings are only unique up to translations on the boundaries.  In fact different choices of gluings lead to nonisometric manifolds (we explain more below).  

It is easier to explain what is going on in the universal covering space.  The universal covering space is built by gluing copies of the $M_{\ga_v}$.  The gluing is described by a tree $\bar {\cal G}$ (the Basse-Serre tree of the fundamental group - see section 2.5 of \cite{C-K1}) with vertex set $\bar {\cal V}$ and edge set $\bar {\cal E}$.  $\Ga$ acts on $\bar {\cal G}$ with $\pi_{\cal G}:\bar {\cal G}\to \bar {\cal G}/\Ga = \cal G$.  We can define the vertex and edge spaces $M_{\bar v}$, $R^2_{\bar e}$ in the universal cover, where $M_{\bar v}$ is isometric to $M_{\ga_v}$ if $\pi_{\cal G}(\bar v)=v$.  For every $s\in \R$ there is an isometry $\tau(\bar v, s):M_{\bar v}\to M_{\bar v}$ which translates by $s$ in the vertical direction.  (The notation $\tau(\bar v,s)$ is somewhat unusual - rather than $\tau^s_{\bar v}$ for example - but it makes things easier to see later.) 

Now let $N$ and $N'$ be nonpositively curved graph manifolds built from the same graph $\cal G$ and isometric pieces $N_v$ and $N_v'$ and let $I_v$ be a chosen isometry respecting the vertical directions.  We can lift the $I_v$'s equivariantly to isometries $I_{\bar v}:M_{\bar v}\to M'_{\bar v}$.  For each $e\in \bar {\cal E}$ from $\bar {v_1}$ to $\bar {v_2}$ then $I_{\bar {v_1}}$ and $I_{\bar {v_2}}$ define two isometries from $R^2_{\bar e}$ to ${R^2}'_{\bar e}$ which differ by a Euclidean translation.  
We can thus uniquely define real numbers $s_{\bar e}$ and $r_{\bar e}$ such that on $R^2_{\bar e}$ we have $$\tau({\bar {v_1}},{s_{\bar e}})\circ I_{\bar {v_1}}=\tau({\bar {v_2}},{r_{\bar e}})\circ I_{\bar {v_2}}.$$
We note that $r_{\bar e}=s_{-\bar e}$ and that the equivariance of the lifting says that if $\pi_{\cal G}(\bar {e_1})=\pi_{\cal G}(\bar {e_2})  $ then $s_{\bar {e_1}}=s_{\bar {e_2}}$ so we can define $s_{e}$ for $e\in \cal E$.  We will call the relative gluings {\it consistent} if for every $v\in {\cal G}$ there is a $s_v\in \R$ such that for every $e\in {\cal E}$ which has $v$ as initial point we have $s_v=s_e$.  If the gluings are consistent then we define an equivariant isometry $I:M\to M'$ by $I(m)=\tau({\bar v},{s_{\bar v}})\circ I_{\bar v}(m)$ for $m\in M_{\bar v}$.  Thus we see:

\begin{lemma}
\label{isometry}
If $N$ and $N'$ are nonpositively curved graph manifolds built from the same graph $\cal G$ and isometric pieces $N_v$ and with consistent relative gluings then $N$ and $N'$ are isometric.
\end{lemma}

\begin{remark}
If the relative gluings are not consistent then it is not hard to see there will be no isometry which takes $N_v$ to $N'_v$.  Thus we can make many nonisometric spaces with the same $\cal G$ and isometric pieces.

Two such nonisometric graph manifolds with the same $\cal G$, and isometric $N_v$, but different gluings will have the same ``geometric data'' in the sense of \cite{C-K1} and hence give a new class of examples of nonisometric graph manifolds whose geometric boundaries are equivariantly homeomorphic. 
\end{remark}

{\bf Proof of theorem \ref{main}}
To prove Theorem \ref{main} we use Remark \ref{key} to define our isometries $I_{\bar v}:M_{\bar v}\to M'_{\bar v}$ and assume (using Lemma \ref{isometry}) that $N$ and $N'$ do not have consistent gluings.  We will then find a closed geodesic $c\in N$ and a curve $\xi \in N'$ representing the same (as identified via the conjugacy $F$) free homotopy class such that $L(\xi)<L(c)$.  This is the desired contradiction since any closed geodesic minimizes the lengths of curves freely homotopic to it and since $\xi$ is freely homotopic to $F(c)$ hence $L(c)=L(F(c))\leq L(\xi)$.  (We will do this by choosing $c$ so that it intersects appropriate boundary tori nearly tangentially and such that the discontinuous image of $c$ when we reglue the Seifert components to get $N'$ can be shortened to get a continuous curve in the appropriate homotopy class.  We actually work in the universal covers.)   

Thus we assume that there are edges $e_1,e_2\in \cal E$ and a vertex $v$ such that $v$ is the final point of $e_1$ and the initial point of $e_2$ and that $s_{-{e_1}}\not = s_{{e_2}}$.  We will let $v_1$ be the initial point of $e_1$ and $v_2$ the final point of $e_2$.

The first step is to find an appropriate closed geodesic $c$ in $N_{v_1}\cup N_{v}\cup N_{v_2}$ such that $c$ intersects each of $T_{e_1}$ and $T_{e_2}$ in exactly two points and intersects no other $T_{e}$.  To choose such a $c$ we first choose closed geodesics $a_i\subset N_{v_i}$ and $b\subset N_v$.  We also choose minimizing geodesics $\tau_i$ from the base point of $a_i$ to the base point of $b$. They are chosen so that $b$ is in the homotopy class of $\gamma_v$ and the classes corresponding to $a_i$ (or more precisely $-\tau_i\cup a_i\cup \tau_i$) do not commute with $\gamma_v$.  For each $n\in \Z$ we consider the curves 
$$c_n=a_1\cup \tau_1 \cup b^n \cup -\tau_2 \cup a_2 \cup \tau_2\cup b^{-n} \cup -\tau_1.$$
We will choose $c$ to be the closed geodesic in the free homotopy class of some $c_n$ (for some $n$ large in absolute value).  Since $c$ can be gotten from $c_n$ by curve shortening, the convexity of the Seifert components can be used to make sure that each curve in the shortening process, and hence the final curve $c$, satisfies the required conditions.

Thus, we can choose a unit speed parameter $c:[0,L]\to N$ so that there are $t_i\in R$, $0=t_0<t_1<t_2<t_3<t_4<t_5=L$ with $c([0,t_1])\subset N_{v_1}$, $c(t_1)\in T_{e_1}$, $c([t_1,t_2])\subset N_{v}$, $c(t_2)\in T_{e_2}$, $c([t_2,t_3])\subset N_{v_2}$, $c(t_3)\in T_{e_2}$, $c([t_3,t_4])\subset N_{v}$, $c(t_4)\in T_{e_1}$, $c([t_4,L])\subset N_{v_1}$.

For such a $c$ we will let $\bar c[0,l]\to M$ be a continuous lift of $c$.  There are thus $\bar v_1,\bar v_2,\bar v_3,\bar v_4,\bar v_5\in \bar {\cal V}$ and $\bar e_1,\bar e_2,\bar e_3,\bar e_4\in \bar {\cal E}$ such that for 
$$i=1,2,3,4,5,\ \ \bar c_i\equiv \bar c([t_{i-1},t_i])\subset M_{\bar v_i},$$ 
$$i=1,2,3,4, \ \  c(t_i)\in \R^2_{\bar e_i},$$ 
$$\pi_{\cal G}(\bar v_1)=\pi_{\cal G}(\bar v_5)=v_1,$$
$$\pi_{\cal G}(\bar v_2)=\pi_{\cal G}(\bar v_4)=v,$$
$$\pi_{\cal G}(\bar v_3)=v_2,$$ 
$$\pi_{\cal G}(\bar e _1)=-\pi_{\cal G}(\bar e_4)=e_1,$$
 and 
$$\pi_{\cal G}(\bar e_2)=-\pi_{\cal G}(\bar e_3)=e_2.$$

We will now determine which $n$ to use.  Let $\R_{\bar v_i}:M_{v_i}\to \R$ represent the projection onto the vertical.
Choose $|n|$ large enough that $|\R_{\bar v_2}(\bar c(t_2))- \R_{\bar v_2}(\bar c(t_1))|>|s_{\bar e_2} - s_{-\bar e_1}|$ (resp. $|\R_{\bar v_4}(\bar c(t_4))- \R_{\bar v_4}(\bar c(t_3))|>|s_{\bar e_2} - s_{-\bar e_1}|=|s_{\bar e_4} - s_{-\bar e_3}|$).  Choose the sign of $n$ (hence the signs of $\R_{\bar v_2}(\bar c(t_2))- \R_{\bar v_2}(\bar c(t_1))$ and $\R_{\bar v_4}(\bar c(t_4))- \R_{\bar v_4}(\bar c(t_3))$ for large enough $n$) so that 
$$|\R_{\bar v_2}(\bar c(t_2))- \R_{\bar v_2}(\bar c(t_1))+ s_{\bar e_2} - s_{-\bar e_1}|<|\R_{\bar v_2}(\bar c(t_2))- \R_{\bar v_2}(\bar c(t_1))|$$

and 

$$|\R_{\bar v_4}(\bar c(t_4))- \R_{\bar v_4}(\bar c(t_3))+ s_{\bar e_4} - s_{-\bar e_3}|<|\R_{\bar v_4}(\bar c(t_4))- \R_{\bar v_4}(\bar c(t_3))|.$$

Having now chosen $c$ we now proceed with the rest of the argument.
We let $\ga \in \Ga$ be such that $\ga(\bar c(0))=\bar c(L)$, and hence in particular 
$$\ga \tau({\bar {v_1}},{s_{\bar e_1}}) (\bar c(0)) = \tau({\bar {v_5}},{s_{\bar e_1}})(\bar c(L)) = \tau({\bar {v_5}},{s_{-\bar e_4}})(\bar c(L))$$.

Consider the following continuous curve, $\xi$ in $M'$ between $\tau({\bar {v_1}},{s_{\bar e_1}}) (\bar c(0))$ and $\tau({\bar {v_5}},{s_{-\bar e_4}})(\bar c(L))$:

For $t_0\leq t\leq t_1$, $\xi(t)=\tau({\bar {v_1}},{s_{\bar e_1}})(\bar c_1(t))$.

For $t_1\leq t\leq t_2$, $\xi(t)=\tau({\bar {v_2}},\frac {t_2-t} {t_2-t_1} {s_{-\bar e_1}}+\frac {t-t_1} {t_2-t_1} {s_{\bar e_2}})(\bar c_2(t))$.

For $t_2\leq t\leq t_3$, $\xi(t)=\tau({\bar {v_3}},{s_{\bar e_3}})(\bar c_3(t))$ (note $s_{\bar e_3}=s_{-\bar e_2}$).

For $t_3\leq t\leq t_4$, $\xi(t)=\tau({\bar {v_4}},\frac {t_2-t} {t_2-t_1} {s_{-\bar e_3}} +\frac {t-t_1} {t_2-t_1} {s_{\bar e_4}}))(\bar c_4(t))$.

For $t_4\leq t\leq t_5$, $\xi(t)=\tau({\bar {v_5}},{s_{-\bar e_4}})(\bar c_5(t))$.

We now only need to see that $L(\xi)<L=L(c)$.  The five pieces $\xi_i$ of $\xi$ correspond to the pieces $c_i$ of $c$.  It is easy to see that $L(c_i)=L(\xi_i)$ for $i=1,3,5$ since one comes from the other by an isometry composed with a vertical shift (also an isometry).  For $i=2$ we note that $\pi_{X_{v_2}}(\bar c_2)$ is the same as $\pi_{X_{v_2}}(\xi_2)$ and hence $L(\xi_2)<L(c_2)$ since:
$$|\R_{\bar v_2}(\xi_2(t_2))- \R_{\bar v_2}(\xi_2(t_1))|=|\R_{\bar v_2}(\bar c_2(t_2))- \R_{\bar v_2}(\bar c_2(t_1))+ s_{\bar e_2} - s_{-\bar e_1}|<$$
$$<|\R_{\bar v_2}(\bar c_2(t_2))- \R_{\bar v_2}(\bar c_2(t_1))|.$$
The $i=4$ case is similar.

\no
Chris Croke:\\
Department of Mathematics\\
University of Pennsylvania\\
209 S. 33rd St.\\
Philadelphia, PA 19104-6395\\
ccroke@math.upenn.edu\\

\end{document}